\documentclass{amsart}
\usepackage{amssymb}
\usepackage{hyperref}
\newcommand{\spli}{{\operatorname{\mathsf   {split}}}}

\newcommand{\suc}{{\operatorname{\mathsf    {succ}}}}
\newcommand{\stem}{{\operatorname{\mathsf   {stem}}}}
\renewcommand{\P}{{\mathbb P}}
\renewcommand{\L}{{\mathbb L}}
\newcommand{\HH}{{\mathbf H}}
\newcommand{\Miller}{{\mathbb M}}

\newcommand{\forces}{\Vdash}
\newcommand{\ZFCa}{{\operatorname{\mathsf {ZFC}}}}
\newcommand{\CH}{\operatorname{\mathsf {CH}}}
\newcommand{\rk}{{{\mathsf {rank}}}}
\newcommand{\reals}{{\mathbb R}}

\newcommand{\rest}{{\mathord{\restriction}}}

\newcommand{\cov}{\operatorname{\mathsf  {cov}}}

\newcommand{\dom}{{\operatorname{\mathsf {dom}}}}

\newcommand{\N}{{\mathcal N}}

\newcommand{\M}{{\mathcal M}}
\newcommand{\V}{{\mathbf V}}
\newcommand{\<}{\langle}
\renewcommand{\>}{\rangle}
\newcommand{\supp}{\mathsf {supp}}
\newcommand{\C}{\mathbf C}
\newcommand{\thinks}{\models}

\newcommand{\lft}[2]{\mathopen\ifcase#1{}\oo\or
                        \big#2\or\Big#2\else\oo\fi} 
\newcommand{\rgt}[2]{\mathclose\ifcase#1{}\oo\or
                        \big#2\or\Big#2\else\oo\fi} 

\newcommand{\SM}{{\mathcal   {SM}}}
\newcommand{\SN}{{\mathcal  {SN}}}

\theoremstyle{plain}
\newtheorem{theorem}{Theorem}
\theoremstyle{plain}
\newtheorem{lemma}[theorem]{Lemma}
\newtheorem{corollary}[theorem]{Corollary}

\newtheorem{definition}[theorem]{Definition}

\begin{document}
\title{Strongly meager sets of size continuum}
\author{Tomek Bartoszynski}
\address{Department of Mathematics and Computer Science\\
Boise State University\\
Boise, Idaho 83725 U.S.A.}
\thanks{First author  partially supported by 
NSF grant DMS 0200671}
\email{tomek@math.boisestate.edu, http://math.boisestate.edu/\char 126 tomek}
\author{Saharon Shelah}
\thanks{Second author partially supported by Israel Science
   Foundation and NSF grant DMS 0072560. Publication 807}
\address{Department of Mathematics\\
Hebrew University\\
Jerusalem, Israel and 
Department of Mathematics\\
Rutgers University\\
New Brunswick, NJ
}
\email{shelah@sunrise.huji.ac.il, http://math.rutgers.edu/\char 126 shelah/}

\begin{abstract}
We will construct several models where there are no strongly meager
sets of size $ 2^{\boldsymbol\aleph_0} $. 
\end{abstract}
\maketitle

\section{Introduction}
In this paper we  work exclusively in the space $2^\omega $
equipped with 
the standard product measure denoted as $\mu$. Let $\N$ and $\M$ denote the
ideal of all $\mu$--measure zero sets, 
and meager subsets of $2^\omega $, respectively.
For $x,y \in 2^\omega$, $x+y \in 2^\omega $ is
defined as $(x+y)(n) = x(n)+y(n) \pmod 2$. In particular, $(2^\omega , 
\operatorname{+})$ is a group and $\mu$ is an invariant measure.

\begin{definition}[\cite{Borel}]
A set $X$ of real numbers or more generally, a metric space, is strong
measure zero  if, for each sequence $\{\varepsilon_n: n \in \omega\}$ of
positive real numbers there is a sequence $\{X_n: n \in \omega\}$ of
subsets of $X$ 
whose union is $X$, and for each $n$ the
diameter of $X_n$ is less than $\varepsilon_n$. 
\end{definition}

The family of strong measure zero
subsets of $2^\omega $ is denoted by $\SN$.

The following characterization of strong measure zero 
 is the starting point for our considerations.
\begin{theorem}[\cite{GMSSmz}]\label{solo}
The following  are equivalent:
\begin{enumerate}
\item $X \in \SN$,
\item for every set $F \in \M$, $X+F \neq 2^\omega$.
\end{enumerate}
\end{theorem}

This theorem indicates that the notion of strong measure
zero should have its category analog.
Indeed, we define after Prikry:
\begin{definition}\label{defstrmea}
  Suppose that $X \subseteq 2^\omega $.
We say that $X$  is strongly meager if for every $H \in
  \N$, $X+H \neq 2^\omega $. Let $\SM$ denote the collection of
  strongly meager sets.
\end{definition}

Observe that if $z \not\in X+F=\{x+f: x \in X, f\in F\}$ then $X \cap
(F+z) = \emptyset$. In particular, a strong measure zero set can be
covered by a translation of any dense $G_\delta $ set, and every
strongly meager set can be covered by a translation of any measure one 
set.

Let Borel Conjecture denote the statement
$\SN=[2^\omega]^{{\mathbf\aleph}_0}$ and Dual Borel Conjecture the
statement: $\SM=[2^\omega]^{{\mathbf\aleph}_0}$.

The question whether Borel Conjecture and Dual Borel Conjecture are
jointly consistent is the motivation for this paper. In particular, we
are interested in what kinds of strongly meager sets exist in various
models.

\section{Preservation of not being strongly meager}
A small modification \cite{kysiak} of the definition of strongly
meager sets leads 
to a family that captures the concept of strongly meager sets while it
has some additional properties.

\begin{definition}
  We say that $X \in \SM^+$ if for every $H \in
  \N$, there exists a countable set $Z \subseteq 2^\omega $ such that
  $$ \forall x \in X \ Z \not \subseteq x+H.$$
\end{definition}

\begin{lemma}
  \begin{enumerate}
  \item $\SM \subseteq \SM^+$,
\item $\SM^+$ is a $\sigma
$-ideal,
\item it is consistent that $\SM=\SM^+$,
\item it is consistent that $\SM\neq\SM^+$.
  \end{enumerate}
\end{lemma}
\begin{proof}
  (1) and (2) are obvious.

(3) $\SM=\SM^+$ holds in Cohen's model, and in fact every known model
where $\SM=[2^\omega]^{{\mathbf\aleph}_0 }$ \cite{CarStr}.

(4) By the results of \cite{BaSh607}, assuming $\CH$, $\SM$ is not a $\sigma$-ideal.
\end{proof}

\begin{definition}
  Suppose that $N \prec \HH(\chi)$ is a countable model. A $G_\delta $
  set $H\in \N$ is big over $N$ if $2^\omega \cap N \subseteq H$.

We say that a proper forcing notion $\P$ is nice if for every $p \in
\P$, a countable model $N \prec\HH(\chi)$ containing $p, \P$ and a set
$H$ which is big over $N$ there exists a condition $q \geq p$ such
that 
\begin{enumerate}
\item $q$ is $(N,\P)$-generic over $N$,
\item $q \forces_\P H \text{ is big over } N[\dot{G}]$.
\end{enumerate}
\end{definition}

\begin{theorem}
  Suppose that $\V \thinks X \not \in \SM^+$ and $\P$ is nice. Then
  $\V^\P \thinks  X \not \in \SM^+$.
\end{theorem}
\begin{proof}
  Suppose that $X\not \in \SM^+$ and let $H \in \N$ be a $G_\delta $
  set such that 
$$\forall Z \in [\V\cap 2^\omega ]^\omega \ \exists x \in X \ Z
\subseteq H+x.$$
Let $\{\dot{x}_n: n \in \omega \}$ be a $\P$-name for a countable set
of reals, and let $p  \in \P$. Find a countable model $N \prec
\HH(\chi)$ containing all relevant objects and let $x \in X$ be such
that $N \cap 2^\omega \subseteq H+x$. Since $\P$ is nice, there exists
$q \geq p$ such that $q \forces_\P \forall n \ \dot{x}_n \in H+x$. 
It follows that $\V^\P \thinks  X \not \in \SM^+$.
\end{proof}

The following theorem is a particular instance of a preservation
result due to Eisworth, \cite{EisDoc}. We give a sketch of the proof
for completeness.
 \begin{theorem}\label{firstpres}
   Let $\{{\mathcal P}_\alpha, \dot{{\mathcal Q}}_\alpha :
\alpha<\delta \}$ be a countable support iteration of proper forcing
notions such that every $\alpha<\delta  $, $\forces_\alpha
\dot{{\mathcal Q}_\alpha} \text{ is nice}$. Then
${\mathcal P}_{\delta}$ is nice.
 \end{theorem}
 \begin{proof}
   The proof follows standard proof of preservation of properness.
\begin{lemma}[\cite{GT}, \cite{BJbook} lemma 6.1.2]\label{newPrelim}
Let 
$\{\mathcal P_\alpha,\dot{\mathcal Q}_\alpha:\alpha<\delta\}$
be a countable support iteration. Assume $\alpha<\beta\leq\delta$,
$p\in\mathcal P_\alpha$ and $\dot\tau$ is a $\mathcal P_\beta$-name for an
ordinal. There is a $\mathcal P_\alpha$-name $\dot\pi$ for an ordinal
and a condition $q\geq_\beta p$ such that 
$q\rest\alpha=p$ and $q\forces_\beta\dot\tau=\dot\pi$.
\end{lemma}

Let $N$ be a countable submodel of $\HH(\chi)$ containing 
$\P_\delta$ and let
$H=\bigcap_n \bigcup_{m>n} U_m$ be big over $N$, where  $U_m$'s are
basic open sets.
\begin{lemma}\label{inducL}
Let
$\{\mathcal P_\alpha,\dot{\mathcal Q}_\alpha:\alpha<\delta\}$
be a countable support iteration and assume that for all
$\alpha<\delta$, 
$\forces_{\alpha}\text{``}\dot{\mathcal Q}_\alpha$ is nice.''
Then for all $\beta\in N\cap\delta$, for all $\alpha\in N\cap\beta$
and for all $p\in\mathcal P_\beta\cap N$, whenever
$q\ge_\alpha p\rest\alpha$ is $(N,\mathcal P_\alpha)$-generic,
there is an $(N,\mathcal P_\beta)$-generic condition 
$r\geq_\alpha p\rest\beta$ such that 
$r\rest\alpha=q$ and $r\forces N[G_\beta] \cap 2^\omega \subseteq H$.
\end{lemma}
\begin{proof} The proof is by induction on $\beta$. 

{\sc Successor step.} Let $\beta=\beta'+1$. We assume that $\beta\in N$
so $\beta'\in N$ too. Using the induction hypothesis on
$\alpha,\beta'$ we can extend $q$ to a 
$(N,\mathcal P_{\beta'})$-generic condition $q'\in\mathcal P_{\beta'}$
such that $q'\ge_{\beta'}p\rest\beta'$
and $q'\rest\alpha=q$. Hence without loss of generality
we can assume $\beta=\alpha+1$.

Since 
$\forces_\alpha N[\dot G_\alpha]\prec \HH(\chi)^{\V[\dot G_\alpha]}$,
$q$ forces that there is an
$(N[\dot G_\alpha],\dot{\mathcal Q}_\alpha)$-generic condition
$\geq\dot p(\alpha)$ which forces that $N[G_\alpha] \cap 2^\omega
\subseteq H$. It has a
$\mathcal P_\alpha$-name, $r(\alpha)$ and
$r=q^\frown r(\alpha)$ is $(N,\mathcal P_\beta)$-generic and $r
\forces N[G_\beta] \cap 2^\omega \subseteq H$.

\vspace{0.1in}

{\sc Limit step.} Let $\beta\in N$ be a limit ordinal and
$\alpha=\alpha_0<\alpha_1<\cdots$ be a sequence of ordinals
in $N$ such that $\sup(N\cap\beta)=\sup_{n\in\omega}\alpha_n$.
Let $\langle\dot\tau_n:n\in\omega\rangle$ enumerate all
$\mathcal P_\beta$-names for ordinals which are in $N$,  let
$\<\dot{x}_n : n \in \omega \> $ enumerate all $\mathcal
P_\beta$-names for reals in $N$ with infinitely many repetitions. 
By induction on $n\in\omega$ using lemma \ref{newPrelim},
we define conditions $p'_n\in\mathcal P_\beta\cap N$, $\mathcal
P_{\alpha_n}$-names $\dot\tau^\star_n\in N$ for  ordinals such
that for each 
$n$,
\begin{enumerate}
\item $p'_0=p$ and $p_{n+1}\geq p_n$,
\item $p'_{n+1}\forces_{\beta}\dot\tau^\star_n=\dot\tau_n$, 
\item $p'_{n+1}\rest\alpha_n=p_n\rest\alpha_n$,
\end{enumerate}

Up to now the proof followed  standard proof of preservation of
properness. Now we will make a small change relevant to the current setup.

We will refine the sequence $\{p_n: n \in \omega\}$ as follows.
Let $\{p^1_m: m \in \omega \}$ and $x^0 \in 2^\omega \cap N$ be such that
\begin{enumerate}
\item $p^1_0=p'_1$,
\item $p^1_{m+1} \geq p^1_m$, $p^1_m \in {\mathcal P}_\delta$,
\item $p^1_m \forces \dot{x}_0 \rest m = x^0 \rest m$.
\end{enumerate}
Since $x^0 \in U_k$ for infinitely many $k$, it follows that 
for some $m$, $p^1_m \forces \dot{x}_0 \in U_m$.  

Let $G$ be a ${\mathcal P}_{\alpha_0}$-generic filter over $N$
containing $p_0 \rest \alpha_0$.
Construct sequence $\{p^2_m: m \in \omega \}$ and $x^1 \in 2^\omega \cap N[G]$
be such that 
\begin{enumerate}
\item $p^2_0= p'_1 \rest [\alpha_0, \delta )$,
\item $p^2_{m+1} \geq p^2_m$, $p^2_m \in {\mathcal P}_{\alpha_0,
    \delta }$,
\item $p^2_m \forces \dot{x}_1 \rest m = x^1 \rest m$.
\end{enumerate}

As before, $N[G] \models x^1 \in U_k $ for infinitely many $k$, so it
follows that 
for some $m > 0$, $N[G] \models p^2_m \forces_{\alpha_0, \delta}
\dot{x_1} \in U_m$.
Thus there is a ${\mathcal P}_{\alpha_0}$-name $\dot{m}$ for an
integer such that 
$$p_2 = p_1 \rest \alpha_0 \star p^2_{\dot{m}} \forces \exists m>0 \ \dot{x}_1 \in U_m.$$
We continue in this fashion and construct a sequence $\<p_n: n \in
\omega \>$ such that
\begin{enumerate}
\item $p_0=p$ and $p_{n+1}\geq p_n$,
\item $p_{n+1}\forces_{\beta}\dot\tau^\star_n=\dot\tau_n$, 
\item $p_{n+1}\rest\alpha_n=p_n\rest\alpha_n$,
\item $p_{n+1} \forces \exists m>n \ \dot{x}_n \in U_m$.
\end{enumerate}  

Note that all conditions $p_n$ are in $N$ (while the sequence is not). 
Now using the induction hypothesis 
 we define
conditions $q_n\in\mathcal P_{\alpha_n}$ such that $q_0=q$,
$q_n\geq p_n\rest\alpha_n$,
$q_{n+1}\rest\alpha_n=q_n$ and
$q_n$ is an $(N,\mathcal P_{\alpha_n})$-generic condition.
We set $r=\bigcup_{n\in\omega}q_n$.
Now for each $n$, since $\supp(p_{n+1})\subseteq N$, 
$\supp(p_{n+1})\subseteq\dom(r)$ and by (1) and (3),
$r\geq p_{n+1}$. Hence $r\forces_\beta\dot\tau_n=\dot\tau^\star_n$,
and 
$q_n=r\rest\alpha_n\forces_{\alpha_n}\dot\tau^\star_n\in N$.
Likewise $r\forces_\beta\dot\tau_n\in N$ for each $n\in\omega$, and
  $r$ is $(N,\mathcal P_\beta)$-generic.
Similarly, for every $n$, $r \forces \exists^\infty m \ \dot{x}_n \in U_m$.
\end{proof}
 \end{proof}

The following definition gives an easy to verify property of forcing
notion which implies niceness.
\begin{definition}
Suppose that $\P $ is a forcing notion satisfying Axiom A.
  $\P $ has $PPP$-property if for every $\P$-name
  $\dot{x}$ such that $\forces_\P \dot{x} \in 2^\omega $ and
  every $p \in \P $, $n \in \omega $ there exists $k \in
  \omega $ such that for every $\ell \in \omega $ there exists $A \in
  [2^\ell]^{<k}$ and $q \geq_n p$, $q \forces_\P
  \dot{x}\rest \ell \in A$. 
\end{definition}

Observe that $PPP$-property is a weak version of Laver property (which
is equivalent to the requirement that $k$ depends only on $n$ but not
on $p$). Thus Laver, Miller and  Sacks forcings have property PPP (since they
have the Laver property), as well as many tree forcings.

\begin{lemma}\label{pppnice}
  If $\P$ has $PPP$-property then $\P$ is nice. In particular, if $\P$
  has Laver property then $\P$ is nice.
\end{lemma}
\begin{proof}
Suppose that $N \prec \HH(\chi)$ is a countable model containing $\P$
and $p \in \P$. Suppose that $H \in \N$ is big over $N$. Let $\<U_j: j
\in \omega\>$ be open sets such that $H=\bigcap_j U_j$.
We want to find a condition $q \geq p$ which is $(N,\P)$-generic such
that $q \forces_\P 2^\omega \cap N[\dot{G}] \subseteq H$.

{\bf Basic step:}
Suppose that $\dot{x}\in N$ is a $\P$-name for a real, $n \in \omega $
and $\bar{p} \in \P \cap N$.
Using $PPP$-property we construct reals $\{x_1, \dots, x_k\}\in N \cap
2^\omega $ and a sequence of conditions $\<\bar{p}_m: m\in \omega
\>\in N$
such that 
\begin{enumerate}
\item $\bar{p}=\bar{p}_0$,
\item $\bar{p}_{m+1} \geq_n \bar{p}_m$,
\item $\bar{p}_m \forces_\P \exists j\leq k \ \dot{x}\rest m = x_j
  \rest m$. 
\end{enumerate}
   
Since sets $U_j$ are open it follows that for every $j$ there is $m$
such that $\bar{p}_m \forces_\P \dot{x} \in U_j$.
In other words, given $\bar{p}$, $n \in \omega $, $j \in \omega $, and $\dot{x}$ as
above we can find $\bar{q} \geq_n \bar{p}$, $\bar{q} \in N$ and
$\bar{q} \forces_\P \dot{x} \in U_j$.

The rest of the construction is standard: we build in $\V$ a sequence
$\<(j_n,\dot{x}_n,A_n): n \in \omega\>$ of all
triples $(j,\dot{x},A)$ where $j \in \omega $, $\dot{x}\in N$ is a
$\P$-name for a real and $A\in N$ is a maximal antichain of $\P$. 
Using the basic step above, we build by induction a sequence $\<p_n: n \in \omega\>$ such
 that 
 \begin{enumerate}
 \item $p_0=p$,
 \item $p_n \in N$,
 \item $p_{n+1} \geq_n p_n$,
 \item $p_{n+1} \forces_\P \dot{x_n} \in U_{j_n}$,
 \item $p_{n+1}$ is compatible with countably many elements of $A_n$.
 \end{enumerate}

Finally, let $q$ be such that $q \geq p_n$ for all $n$. Clearly, $q$
has the required properties.
\end{proof}

\section{No small sets of size continuum}
Suppose that ${\mathcal J}$ is a $\sigma$-ideal of sets of
reals having Borel basis. Define
${\mathcal J}^\star = \{X \subseteq 2^\omega: \forall A \in {\mathcal
  J} \ X+A \neq 2^\omega \}$.
In particular, $\M^\star = \SN$ and $\N^\star=\SM$.

 \begin{definition}
   We say that ${\mathcal P}$ has property $P_{\mathcal J}$
 if for every condition $p \in {\mathcal P}$, and a
  ${\mathcal P}$-name for a real $\dot{x}$ such that $p
  \forces_{\mathcal P}  \dot{x} \not\in \V\cap 2^\omega $ there exists a
  set $H \in \V \cap {\mathcal J}  $ such that for all $z \in \V \cap
  2^\omega$, there exists $q_z \geq p$ such that 
$$q_z\forces \dot{x} \in H+z.$$
 \end{definition}

 \begin{theorem}\label{basic}
Let $\{{\mathcal P}_\alpha, \dot{{\mathcal Q}}_\alpha :
\alpha<\omega_2\}$ be a countable support iteration of proper forcing
notions.
If for every $\alpha<\omega_2 $, ${\mathcal
     P}_{\omega_2}/{\mathcal P}_\alpha $ has property $P_{\mathcal J}$, then
$\V^{{\mathcal P}_{\omega_2}} \thinks {\mathcal J}^\star \subseteq
[\reals]^{<2^{\boldsymbol\aleph_0}} .$
 \end{theorem}
 \begin{proof}
   Suppose that $X \subseteq 2^\omega \cap \V^{{\mathcal
       P}_{\omega_2}}$ and $\forces_{\omega_2} X \in {\mathcal J}^\star$.
Working in $\V^{{\mathcal P}_{\omega_2}}$, for every Borel set $H \in
{\mathcal J}  $ find a real $x_H$ such that $\V^{{\mathcal P}_{\omega_2}}
\thinks x_H \not\in H+X$.
Granted that for every $\beta<\omega_2$, $\V^{{\mathcal P}_\beta}
\thinks \CH$, we can find $\alpha <\omega_2$ such that $x_H \in
\V^{{\mathcal P}_\alpha}$ for $H \in \V^{{\mathcal P}_\alpha }$.
Suppose that $p \forces_{\omega_2} \dot{x} \not\in \V^{{\mathcal
    P}_\alpha} \cap 2^\omega $.
Apply property $P_{\mathcal J}$ to find $H \in {\mathcal J} $ with the required properties.
In particular, there is $q \geq p$ such that $q \forces_{\omega_2}
\dot{x} \in H+x_H$. It follows that $q \forces_{\omega_2} x_H \in H+
\dot{x}$, thus $q \forces_{\omega_2} \dot{x} \not\in X$. 
Since $\dot{x}$ and $p$ were arbitrary, we conclude that $X \subseteq
\V^{{\mathcal P}_\alpha} \cap 2^\omega $.
 \end{proof}

 \begin{definition}
   Suppose that $\P$ satisfies Axiom A. We say that $\P$ is strongly
   $\omega^\omega $-bounding if for every $p \in \P$, $n \in \omega $
   and a $\P$-name $\dot{n}$, if $p \forces \dot{n} \in \omega $ then
   there exists $q \geq_n p$ and a set $A \in [\omega]^{<\omega}$ such
   that  $q \forces \dot{n} \in A$.
 \end{definition}

The following is (essentially) proved in \cite{GJS92}, see also
\cite{BJbook} theorem 8.2.14. 
 \begin{theorem}
   If $\P$ is strongly $\omega^\omega $-bounding then $\P$ has
   property $P_\M$. Let $\{{\mathcal P}_\alpha, \dot{{\mathcal Q}}_\alpha :
\alpha<\omega_2\}$ be a countable support iteration of proper forcing
notions such that  for $\alpha < \omega_2$, $\forces_\alpha \dot{{{\mathcal Q}}}_\alpha$ is
   strongly $\omega^\omega $-bounding. 
Then
$\V^{{\mathcal P}_{\omega_2}} \thinks \M^\star=\SN \subseteq
[\reals]^{<2^{\boldsymbol\aleph_0}} .$
 \end{theorem}

Laver forcing (and many other forcing notions) have property $\P_\N$,
a somewhat weaker result was proved in \cite{NowWei},
but whether this property is preserved, particularly for posets adding
unbounded reals, is unclear. 
Therefore we will use the following notions.
\begin{definition}
   Suppose that $N \prec \HH(\chi)$. We say that a sequence of clopen
   sets $\C=\<C_n: n \in \omega\>$ is big over $N$ if
\begin{enumerate}
\item $C_n$'s have pairwise disjoint supports,
\item $\mu(C_n)\leq 2^{-n}$ for $n \in \omega $,
\item  for every infinite set $X \subseteq 2^\omega$, $X \in N$, 
there exists infinitely many $n$ such that  $ X+C_n=2^\omega $.
\end{enumerate}

 \end{definition}

 \begin{lemma}
   Suppose that $\<C_n: n \in \omega\>$ is big over $N$. Then for
   every $x \in 2^\omega $,
   $\<C_n+x: n \in \omega\>$ is big over $N$.
 \end{lemma}
 \begin{proof}
For a clopen set $C$,
   if $C+X=2^\omega $ then $C+x+X=2^\omega +x=2^\omega $.
 \end{proof}
 \begin{theorem}
   For every countable model $N \prec \HH(\chi)$ there exists a set $H
   \in \N$ which is big over $N$.
 \end{theorem}
 \begin{proof}
We will need the following theorem:
\begin{theorem}[\cite{ErdKunMaul81Som}, \cite{BJbook} theorem 3.3.9]\label{lorenz}
Suppose that $X \subseteq 2^\omega$ is an infinite set and
$\varepsilon>0$. Then there exists a clopen set $C \subseteq 2^\omega$
such that $\mu(C)<\varepsilon$ and $X+C=2^\omega$.
\end{theorem}

Let $N \prec \HH(\chi)$ be a countable model.
Fix an enumeration (with infinitely many repetitions), $\{X_n: n \in
\omega \}$ of all infinite subsets of 
$N\cap 2^\omega $ that are in $N$.
By theorem \ref{lorenz}, for each $ n$ there exists a clopen set
$C_n$ of measure $<2^{-n}$ such that $C_n + X_n=2^\omega$. 
 \end{proof}

 \begin{corollary}
   Suppose that $c$ is a Cohen real over $\V$. There exists a sequence
   $\<C_n: n \in \omega\> \in \V[c]$ which is big over $\V$. 
 \end{corollary}

 \begin{definition}
We say that a proper forcing notion $\P$ is good if for every $p \in
\P$, a countable model $N \prec\HH(\chi)$ containing $p, \P$ and a
sequence $\C=\<C_n: n \in \omega\>$
which is big over $N$ there exists a condition $q \geq p$ such
that 
\begin{enumerate}
\item $q$ is $(N,\P)$-generic over $N$,
\item $q \forces_\P \<C_n:n \in \omega\> \text{ is big over }
  N[\dot{G}]$,
\item $q\forces_\P F_\C\cap N=F_\C \cap N[\dot{G}]$, where
$F_\C=\{x\in 2^\omega : \forall^\infty n \ x \not\in C_n\}$.
\end{enumerate}

 \end{definition}

 \begin{theorem}\label{21}
   If $\{{\mathcal P}_\alpha, \dot{{\mathcal Q}}_\alpha:
   \alpha<\omega_2\}$ is a countable support iteration and for every
   $\alpha<\omega_2$, $\forces_\alpha \text{``} \dot{{\mathcal
       Q}}_\alpha $ is good'' then ${\mathcal P}_{\omega_2}$ is good.
 \end{theorem}
 \begin{proof}
   Similar to the proof of theorem \ref{firstpres}. Note that the
   requirements are represented by $G_\delta $ sets. 
 \end{proof}

 \begin{theorem}
   Let $\{{\mathcal P}_\alpha, \dot{{\mathcal Q}}_\alpha :
\alpha<\omega_2\}$ be a countable support iteration of proper forcing
notions such that every $\alpha<\omega_2 $, $\forces_\alpha
\dot{{\mathcal Q}}_\alpha$ is good. Then
$\V^{{\mathcal P}_{\omega_2}} \thinks \SM \subseteq
[\reals]^{<2^{\boldsymbol\aleph_0}} .$
 \end{theorem}
 \begin{proof}
   Suppose that $X \subseteq 2^\omega \cap \V^{{\mathcal
       P}_{\omega_2}}$ and $\V^{{\mathcal P}_{\omega_2}} \thinks X \in
   \SM$.
As in the proof of theorem \ref{basic}, for $H \in \N$ let $x_H \in
2^\omega $ be such that $x_H \not\in H+X$.
We can assume that $x_H \in \V$ if $H \in \V$.
Suppose that $\dot{x}$ is a ${\mathcal P}_{\omega_2}$-name for a real
and $p \forces_{{\mathcal P}_{\omega_2}} \dot{x}\not\in \V$.
Let $N \prec\HH(\chi)$ be a countable model containing $p, {\mathcal
  P}_{\omega_2}$, and $\dot{x}$.
Let $\<C_n:n \in \omega\> \in \V$ be big over $N$.
Put $H=\bigcap_m \bigcup_{n>m} C_n \in \N$ and let
$\C=\<C_n+x_H: n \in \omega\>$. Sequence $\C$ is also big over $N$,
and $H+x_H = \bigcap_m \bigcup_{n>m} C_n+x_H$.
Since ${\mathcal P}_{\omega_2}$ is good we can find $q \geq p$ such
that 
\begin{enumerate}
\item $q$ is $(N, {\mathcal P}_{\omega_2})$-generic over $N$,
\item $q \forces_{\omega_2} \C \text{ is big over } N[\dot{G}]$,
\item $q \forces_{\omega_2} F_\C \cap N=F_\C \cap N[\dot{G}]$.
\end{enumerate}

Since $p \forces \dot{x} \not\in N$ it follows that $q \forces \dot{x}
\not \in F_\C=2^\omega \setminus (H+x_H)$. In particular, $q \forces
\dot{x} \in H+x_H$, and as before, $q \forces \dot{x} \not\in X$.
Thus $X \subseteq \V \cap 2^\omega $, which finishes the proof.
 \end{proof}

\section{Laver and Miller forcings}
In this section we will show that Laver and Miller forcings are good,
and therefore there are no strongly meager sets of size continuum in
Laver's or Miller's model.

For a tree $p \subseteq 2^{<\omega}$ let $[p]$ denote the set of
branches of $p$, and let $p_s=\{t \in p: s \subseteq t \text{ or } t
\subseteq s\}$. 
For $s \in p$ and $m>|s|$ let $\suc^m_p(s)=\{t \in p : |s|+|t|=m
\text{ and }
s^\frown t \in p\}$. We will suppress the superscript $m$ if its value
will be clear from the context. 

Let $\spli(p) = \{s\in p : |\suc_p(s)|>1\} = \bigcup_{n \in \omega}
\spli_n(p),$ where
$\spli_n(p)=\left\{s \in \spli(p): \lft2|\lft1\{t \subsetneq s: t \in
\spli(p)\rgt1\}\rgt2| =n\right\}$. The unique element of $\spli_0(p)$
is called $\stem(p)$.  For $s \in p$ let
$A_s = \{n \in \omega : s^\frown n \in p\}$

\begin{definition}
The Laver forcing $\L$ is the following forcing notion:
$$p \in \L \iff p \subseteq \omega^{<\omega} \text{ is a
  tree } \& \ 
\forall s \in p \ \lft1(|s| \geq \stem(p) \rightarrow
|\suc_p(s)|=\boldsymbol\aleph_{0}\rgt1).$$

For $p, q \in \L$, $p \geq q$ if $p \subseteq q$.

For every $p \in \L$ and $s \in \omega^{<\omega}$ define a node
$p(s)$ in the following way:
$p(\emptyset)=\stem(p)$ and for $n \in \omega$ 
let $p(s^\frown n)$ be the $n$-th 
element of $\suc_p\lft1(p(s)\rgt1)$.

For $p,q \in \L$  and $n \in \omega$ define
$$p \geq_n q \iff p \geq q \ \&\ 
\forall s \in n^{\leq n}\ 
p(s)=q(s).$$
In particular, $p \geq_0 q$ is equivalent to $p \geq q$ and
$\stem(p)=\stem(q)$.

The rational perfect forcing (Miller forcing) 
$\Miller$ is the following forcing notion:
\begin{multline*}
p \in \Miller \iff p \subseteq \omega^{<\omega} \text{ is a
  perfect tree } \& \\ 
\forall s \in p \ \exists t \in p \
\lft1(s \subseteq t \ \& \ |\suc_p(t)|=\boldsymbol\aleph_{0}\rgt1).
\end{multline*} 
For $p, q \in \Miller$, $p \geq q$ if $p \subseteq q$.
Without loss of generality we can assume that  $|\suc_p(s)|=1$ or
$|\suc_p(s)|=\boldsymbol\aleph_0$ for all $p\in \Miller$ and $s \in
p$.
Conditions of this type form a dense subset of $\Miller$.

For $p,q\in \Miller$, $n \in \omega$,   we let 
$$p \geq_n q \iff p \geq q \ \& \ \spli_n(q)=\spli_n(p).$$

\end{definition}

 \begin{theorem}
Laver and Miller forcings are good.
 \end{theorem}
 \begin{proof}
We will prove the theorem for the Laver forcing $\L$, the proof for
Miller forcing is similar (and easier).

\begin{lemma}\label{lem1}
Suppose that $p \forces_\L \dot{x} \in 2^\omega $.
  There exists $q \geq_0 p$ a set of reals $\{x_s : s \in q,\ 
\stem(q) \subseteq s\}$ such that  
\begin{enumerate}
\item if $\stem(q) \subseteq s \in q$ then for every $n$ there is
  $r\geq_0 q_s$ such that  $r \forces_\L \dot{x} \rest n = x_s \rest
  n$,
\item for every $s$, mapping $n \leadsto x_{s^\frown n}$ ($n \in A_s$) is either one-to-one or
  constant,
\item $\lim_{n \in A_s} x_{s^\frown n} =x_s$.
\end{enumerate}
\end{lemma}
\begin{proof}
  Induction on levels of $p$. Suppose that $|\stem(p)|=k$ and define
  trees $\{p_n: n \in \omega \}$ such that  
  \begin{enumerate}
  \item $p_0=p$,
  \item $p_{n+1} \cap \omega^{k+n} = p_n \cap \omega^{k+n}$ for every
    $n$,
  \item reals $x_s$ are defined for $s \in  p_{n+1} \cap \omega^{k+n}$
    and satisfy conditions (1) and (3).
  \end{enumerate}
Suppose that $p_n$ is given, and $s \in p_n \cap \omega^{k+n}$. 
For each $l \in A_s$ let $r_l^s \geq_0 (p_n)_{s^\frown n}$ and $s_l \in
2^l$ be such that  $r_l^s \forces \dot{x} \rest l = s_l$.
Let $A'_s \subseteq A_s$ be an infinite set such that  $\bigcup_{l \in
  A'_s} s_l = x_s \in 2^\omega $. Finally let 
$p_{n+1} = \bigcup_{s \in p_n \cap \omega^{k+n}} \bigcup_{l \in A'_s}
r_l^s$.
Clearly, $q'=\bigcap_n p_n$ satisfies (1) and (3). By pruning it
further we get $q \geq_0 q'$ that satisfies (2) as well.
\end{proof}

Recall that $B \subseteq p$ is a front in $p$, that is for every
branch $x \in [p]$ 
there is $n \in \omega $ such that  $x \rest n \in B$. 
\begin{lemma}\label{lem2}
  Suppose that $p \in \L$ and $B \subseteq p$. Then either there
  exists $q \geq_0 p $ such that  $q \cap B=\emptyset$ or there exists
  $q \geq_0 p$ such that  $B$ is a front in $q$.
\end{lemma}
\begin{proof}
  Define rank function on $p$ as follows, for $s \in p$
  \begin{enumerate}
  \item $\rk_B(s)=0$ if $\exists n \leq |s|\ s \rest n \in B$,
  \item If $\rk_B(s)\neq 0$ then $$\rk_B(s)=\min\{\alpha : \exists A
    \in [A_s]^\omega \ \forall n \in A \ \rk_B(s^\frown n)<\alpha \}.$$
  \end{enumerate}
If $\rk_B(\stem(p))$ is defined then let $q \geq_0 p$ be such that
$\rk_B$ is defined on $q$ and $\rk_B(s) > \rk_B(t)$ whenever $t
\subsetneq s$. It follows that $B$ is a front in $q$.

 If $\rk_B(\stem(p))$ is not defined then there is $q \geq_0 p$ such
 that  $\rk_B(s)$ is undefined for $s \in q$. Thus $q \cap B = \emptyset$.
\end{proof}

   Suppose that $N \prec \HH(\chi)$,  $\<C_n:n
   \in \omega\>$ is big over $N$.

We will need the following two basic constructions.

{\bf Basic Step I}. Given $p \in N \cap \L$, $k,m \in \omega $ 
and an $\L$-name for a
real $\dot{x} \in N$ such
that  $p \forces_\L \dot{x} \not\in N \cap 2^\omega $ we need $q
\geq_k p$ and $n > m$ such that  $q \forces_\L \dot{x} \in C_n$.

We can
assume that $k=0$, otherwise we repeat the construction described
below for each ``protected'' node.

Working in $N$ we can shrink $p$ (without changing the stem) so that
it satisfies conclusion of lemma
\ref{lem1}.
Let $B=\{s \in p: \forall i,j \in A_s \ x_{s^\frown i}\neq x_{s^\frown
  j}\}$.  Applying lemma \ref{lem2} we can further shrink $p$ (without
changing the stem) so that $B$ is a front in $[p]$. Note that if 
$p \cap B = \emptyset$ then $p \forces_\L \dot{x} = x_{\stem(p)}$,
which contradicts the choice of $\dot{x}$.

Let $Z \subseteq B$ be a maximal antichain in $p$. Fix $s \in Z$ and let $X_s
= \{x_{s^\frown l}: l \in A_s\}$. Since $N \models X_s \in [2^\omega
]^\omega $ the set
$$B_s=\{l \in A_s: \exists n > m \ x_{s^\frown l} \in C_n\}$$
is infinite. 
For each $l \in B_s$ fix $n_l$ such that  $x_{s^\frown l} \in
C_{n_l}$. Let $q^s_l \geq_0 p_{s^\frown l}$ be such that  $q^s_l
\forces_\L \dot{x} \rest \supp(C_{n_l})=x_{s^\frown l}\rest
\supp(C_{n_l})$. In particular, $q^s_l
\forces_\L \dot{x} \in C_{n_l}$. 

Let $q_s = \bigcup_{l \in B_s} q^s_l$ and let 
$q = \bigcup_{s \in Z} q_s$. Clearly $q \geq_0 p$ and 
$q \forces_\L \exists n>m \ \dot{x} \in C_n$.

{\bf Basic Step II}. Given $p \in N \cap \L$, $k,l \in \omega $ 
and an $\L$-name $\dot{X}=\{\dot{x}_j: j \in \omega\}$ for a countable
set of reals in $N^{\L}$ we need $q \geq_k p$ and $n>l$ such that  
$q \forces_\L \dot{X}+C_n=2^\omega $.

As before we can assume that $k=0$. Define
sequences $\<q_j: j \in \omega \>$, $\{s^j_i: i \leq j, j \in \omega \}$ 
such that  
\begin{enumerate}
\item $q_0=p$,
\item $q_{j+1} \geq_0 q_j$,
\item $|s^j_i| > j$ for $i \leq j$,
\item $s^j_u\neq s^j_v$ if $u,v \leq j$, $u \neq v$,
\item $s^j_u \subseteq s^i_u$ if $u \leq j \leq i$.
\item $q_{j+1} \forces_\L
  \forall i \leq n \ s^j_i \subseteq \dot{x_i}$,
\end{enumerate}

Let $x_i = \bigcup_{j>i} s^j_i$. Reals $x_i$ are pairwise different
and $X=\{x_i : i \in \omega\} \in N$. By the assumption, there exists
$n>l$ such that  $X+C_n=2^\omega$. Find $l \in \omega $ such that
$C_n \subseteq 2^l$ and let $j$ be so large that $q_j \forces_\L \{x
\rest l : x\in X\} = \{s^j_i \rest l: i \leq j\}$. Clearly,
$q_j \forces_\L \dot{X}+C_n=2^\omega $.

The rest of the proof is standard induction. We proceed as in the
proof of lemma \ref{pppnice}. Let $\{(\dot{X}_n, \dot{x}_n, A_n): n \in
\omega \}$ be an enumeration of all sets $(\dot{X},\dot{x}, A)$ in $N$ where
$\dot{X}$ is an $\L$-name for a countable sets of reals in $N$,
$\dot{x}$ is a $\L$-name for 
a real such that  $\forces \dot{x} \not\in N \cap 2^\omega $, $A$ is a
maximal antichain in $\L$.
Using the basic steps above, we build by induction a sequence $\<p_n:
n \in \omega\>$ such 
 that 
 \begin{enumerate}
 \item $p_0=p$,
 \item $p_{n+1} \geq_n p_n$,
 \item $p_{n+1} \forces \exists m>n \ \dot{x_n} \in C_m$,
 \item $p_{n+1} \forces \exists m>n \ \dot{X}_n+ C_m=2^\omega $,
 \item $p_{n+1}$ is compatible with countably many elements of $A_n$.
 \end{enumerate}
 Note that conditions produced in step I are not necessarily in $N$, but
 they are of form $\bigcup_{s \in Z} q_s$ where $Z \in N$ and
 $q_{s^\frown l} \in
 N$ for $s \in Z$, $l \in \omega $. This suffices to carry out the
 construction. 
 \end{proof}

It is not clear whether every forcing notion having Laver property (or
PPP property) is nice, but many of them are.

\section{Strongly meager sets of size continuum}
The simplest case when we have large strongly meager sets is when
there are many random reals around.
\begin{lemma}
  If $\cov(\N)=2^{{\mathbf \aleph}_0}$ then $[2^\omega ] \subsetneq \SM$.
\end{lemma}
\begin{proof}
  Suppose that $|X|<2^{{\mathbf \aleph}_0}$ and $H \in \N$. It follows that
$H+X = \bigcup_{x\in X} H+x \neq 2^\omega$. Therefore $X \in \SM$.

To construct a strongly meager set of size continuum fix an
enumeration $\{H_\alpha : \alpha < 2^{{\mathbf \aleph}_0}\}$ of Borel
measure zero sets and construct
inductively reals $\{x_\alpha, z_\alpha : \alpha < 2^{{\mathbf
    \aleph}_0}\}$ such that  
\begin{enumerate}
\item $\{x_\beta: \beta<\alpha \} \cap (H_\alpha +z_\alpha ) =
  \emptyset$,
\item $x_\alpha \not \in \bigcup_{\beta<\alpha } H_\beta+z_\beta$. 
\end{enumerate}
Note that $X=\{x_\alpha : \alpha < 2^{{\mathbf \aleph}_0}\}$ is the
required set since $z_\beta \not\in H_\beta +X$.
\end{proof}

\begin{theorem}
  It is consistent with $\ZFCa$ that $\cov(\N)< 2^{{\mathbf
      \aleph}_0}$ and there exists a strongly meager set of size continuum.
\end{theorem}
\begin{proof}
  Let $\V'$ be a model obtained by adding ${\mathbf \aleph}_2$ Cohen
  reals to a model $\V_0 \models \CH$, and
  let $\V$ be obtained by adding ${\mathbf \aleph}_1$ random
  reals (side-by-side) to $\V'$.
It is easy to see that $\V \models \V' \cap 2^\omega \in \SM$ (as
witnessed by random reals), and it
is well known (\cite{Paw86Sol}, \cite{BJbook} theorem 3.3.24) that $\V \models \cov(\N)={\mathbf \aleph}_1$.
\end{proof}


\end{document}